\documentclass[12pt,leqno]{article}

\usepackage{amsmath,amsthm,amsfonts,latexsym,amscd,amssymb}
\numberwithin{equation}{section}
\input xypic
\def\qed{{\hbadness=10000\hfill\ \vbox{\hrule height.09ex
   \hbox{\vrule width.09ex height1.55ex depth.2ex \kern1.8ex
   \vrule width.09ex height1.55ex depth.2ex}\hrule height.09ex}\break
   \bigskip}}
\newtheorem{theorem}{Theorem}[section]
\newtheorem{lemma}{Lemma}[section]
\newtheorem{corollary}{Corollary}[section]

\theoremstyle{definition}

\theoremstyle{remark}

\begin{document}

\linespread{1}\title{\textbf{Hypersurface of a Finsler space subjected to an $\textsl{h}$-exponential change of metric}}

\author{M.$\,$K. \textsc{Gupta}\thanks{Supported by UGC, Government of India}\, and Anil K.\,\textsc{Gupta}\\
\normalsize{Department of Pure $\&$ Applied Mathematics}\\
\normalsize{Guru Ghasidas Vishwavidyalaya}\\
\normalsize{Bilaspur (C.G.), India}\\
\normalsize{Email: mkgiaps@gmail.com\,; gupta.anil409@gmail.com}}
\date{}
\maketitle

\linespread{1.3}\begin{abstract}
Recently we have obtained the Cartan connection for the Finsler space whose metric is given by an exponential change with an \textsl{h}-vector. In this paper, we discuss certain geometric properties of a Finslerian hyperspace subjected to an $\textsl{h}$-exponential change of metric.\\   
\textbf{Keywords}\,:\,Finsler space, hypersurface, exponential change, $\textsl{h}$-vector.\\
2000 Mathematics Subject Classification\,:\textbf{ 53B40}.
\end{abstract}
\section{Introduction}
~~~~In 2006, YU Yao-yong and YOU Ying \cite{yy06} studied a Finsler space with metric function given by exponential change of Riemannian metric. In 2012, H.\,S.\,Shukla et.\,al. \cite{hb12} considered a Finsler space $ \overline F^n=(M^n,\overline L)$, whose Fundamental metric function is an exponential change of Finsler metric function given by \\[-8mm] 
\begin{equation*}
\overline L=L\,e^{\beta /L},
\end{equation*}\\[-10mm]
 where $\beta=b_i(x)y^i$ is $1$-form on manifold $M^n$.

H.\,Izumi \cite{hi80} introduced the concept of an \textsl{h}-vector $b_i(x,y)$ which is \textsl{v}-covariant constant with respect to the Cartan connection and satisfies $L\,C^h_{ij}\,b_h=\rho\,h_{ij}$, where $\rho$ is  a non-zero scalar function and $C^i_{jk}$ are components of Cartan tensor.\,Thus if $b_i$ is an \textsl{h}-vector, then\\[-12mm]
\begin{equation}
(i)\,\, b_i|_k=0,\quad \quad\quad (ii)\,\,L\,C^h_{ij}b_h=\rho\, h_{ij}\,.
\end{equation}
From the above definition, we have \\[-12mm]
\begin{equation} L\,\dot\partial_j b_i= \rho\, h_{ij}\,,\end{equation}\\[-12mm]
which shows that $b_i$ is a function of directional argument also. H.\,Izumi \cite{hi80} proved that  the scalar $\rho$ is independent of directional argument.\,Gupta and Pandey \cite{mp15} proved that if the \textsl{h}-vector $b_i$ is gradient then the scalar $\rho$ is constant

B.$\,$N.$\,$Prasad \cite{bn90} obtained the Cartan connection of Finsler space whose metric is given by \textsl{h}-Rander's change of a Finsler metric.\,Gupta and Pandey \cite{mp15} obtained the Cartan connection of Finsler space whose metric is given by \textsl{h}-Kropina change of Finsler metric. Present authors \cite{ma16} studied the Cartan connection of Finsler space whose metric is given by \textsl{h}-exponential change of Finsler metric.

The theory of hypersurfaces in a Finsler space has been introduced by by E. Cartan \cite{ec34}. A. Rapcs\`ak \cite{ar57} introduced three kinds of hyperplanes and M.\,Matsumoto \cite{mm85} has classified the hyersurfaces and developed a systematic theory of Finslerian hypersurfaces. Gupta and Pandey \cite{mp08,mp09} discussed the hypersurface of a Finsler space whose metric is given by certain transformation with an \textsl{h}-vector.
 
In the present paper, we discuss the geometric properties of hypersurface of a Finsler space $^*\!F^n=(M^n,{^*}\!L)$, whose metric function  $^*\!L$\,is given by an \textsl{h}-exponential change of a Finsler metric function \textit{i.e.}\,\\[-10mm]
\begin{equation}
^*\!L=L\,e^\frac{\beta}{L}\,,
\end{equation}\\[-10mm]
where $\beta=b_i(x,y)y^i $ and $b_i$ is an \textsl{h}-vector.
\section{Preliminaries}

~~~~Let $F^n=(M^n,L)$ be an $n$-dimensional Finsler space equipped with the Fundamental function $L(x,y)$. The metric tensor, angular metric tensor and Cartan tensor are defined by $g^{}_{ij}=\frac{1}{2}\dot\partial_i\dot\partial_jL^2$, $h_{ij}=g^{}_{ij}-l_il_j$ \,and\, $C_{ijk}=\frac{1}{2}\dot\partial_ig_{jk}$ respectively, where $\dot\partial_k=\frac{\partial}{\partial y^k}$\,.\,The Cartan connection is given by $C\Gamma = (F^i_{jk},N^i_k,C^i_{jk})$. The \textsl{h}- and \textsl{v}-covariant derivatives $X_{i|j}$ and $X_i|_j$ of a covarient vector field $X_i$ are defined by 
\begin{equation} X_{i|j} =\partial_j X_i -N^r_j\,\dot\partial_r X_i-X_r F^r_{ij}\,,\end{equation}\\[-12mm]
and \\[-12mm]
\begin{equation} X_i|_j= \dot\partial_j X_i -X_rC^r_{ij}\,,\end{equation}
where $\partial_k=\frac{\partial}{\partial x^k}$\,.

A hypersurface $M^{n-1}$ of the underlying smooth manifold $M^n$ may be parametrically represented by the equation $x^i = x^i(u^\alpha)$, where $u^\alpha$ are Gaussian coordinates on $M^{n-1}$ (Latin indices run from 1 to $n$ while Greek indices run from 1 to $n-1$). Here, we shall assume that the matrix consisting of the pojection factors $B^i_\alpha = \partial x^i/\partial u^\alpha$ is of rank $n-1$. If the supporting element $y^i$ at a point $u=(u^\alpha)$ of $M^{n-1}$ is assumed to be tangent to $M^{n-1}$, we may then write $y^i=B^i_\alpha(u)v^\alpha$ so that $v=(v^\alpha)$ is thought of as the supporting element of $M^{n-1}$ at a point $u^\alpha$. Since the function $\underline L{(u,v)}=L(x(u),y(u,v))$ gives arise a Finsler function on $M^{n-1}$, we get an ($n$-1)-dimensional Finsler space $F^{n-1}=(M^{n-1},\underline L(u,v))$. 

At each point $u^\alpha$ of $F^{n-1}$, the unit normal vector $N^i(u,v)$ is defined as
\begin{equation}
g_{ij}B^i_\alpha N^j=0 ~~~~~~~ g_{ij}N^iN^j=1\,,
\end{equation}
The inverse projection factors $B^\alpha_i(u,v)$ of $B^i_\alpha$ are defined as 
\begin{equation}
B^\alpha_i=g^{\alpha\beta}g_{ij}B^j_\alpha\,,
\end{equation}
where $g^{\alpha\beta}$ is the inverse of metric tensor $g_{\alpha\beta}$ of $F^{n-1}$.\\
from (2.3) and (2.4), it follows that
\begin{equation}
B^i_\alpha B^\beta_i=\delta^\beta_\alpha\,, ~~~~B^i_\alpha N_i=0\,, ~~~~~N_iB^i_\alpha=0\,,~~~~~ N_i N^i=1,
\end{equation}\\[-12mm]
and further\\[-12mm]
\begin{equation}
B^i_\alpha B^\beta_j+ N^i N_j=\delta^i_j\,.
\end{equation}
For the induced Cartan connection $ ICT=(F^\alpha_{\beta \gamma},G^\alpha_\beta,C^\alpha_{\beta \gamma})$ on $F^{n-1}$, the second fundamental $h$-tensor $H_{\alpha \beta}$ and the normal curvature vector $H_\alpha$ are given by
\begin{equation}
H_{\alpha \beta} =N_i(B^i_{\alpha \beta}+F^i_{jk}B^j_\alpha B^k_\beta)+M_\alpha H_\beta \,,
\end{equation}\\[-12mm]
and\\[-12mm]
\begin{equation}
H_\alpha=N_i(B^i_{0\alpha}+G^i_j B^j_\alpha)\,,
\end{equation}
where $M_\alpha=C_{ijk}B^i_\alpha N^j N^k$,~~$B^i_{\alpha\beta}=\frac{\partial^2x^i}{\partial u^\alpha \partial v^\beta}$ and $B^i_{0\alpha}=B^i_{\beta \alpha}v^\beta$.\\
The equations (2.7) and (2.8) yield
\begin{equation}
H_{0\alpha}=H_{\beta\alpha}v^\beta =H_\alpha,~~~~~H_{\alpha 0}=H_{\alpha\beta}v^\beta =H_\alpha +M_\alpha H_0\,.
\end{equation}\\[-10mm]
The second fundamental \textsl{v}-tensor $M_{\alpha\beta}$ is defined as\\[-10mm]
\begin{equation}
M_{\alpha\beta}=C_{ijk} B^i_\alpha B^j_\beta N^k \,.
\end{equation}
The relative \textsl{h}- and \textsl{v}-covariant derivatives of $B^i_\alpha$ and $N^i$ are given by
\begin{equation}\begin{split}
& B^i_{\alpha|\beta}=H_{\alpha\beta} N^i\,,~~~~ ~~~~~ B^i_\alpha|^{}_\beta=M_{\alpha\beta} N^i\,, ~~~~~~~\\& N^i_{{}_{|\beta}}=-H_{\alpha\beta}B^\alpha_j g^{ij}\,, ~~~~~~~ N^i|^{}_\beta=-M_{\alpha\beta}B^\alpha_j g^{ij}\,.
\end{split}\end{equation}
Let $X_i(x,y)$ be a vector field on $F^n$.\,Then the relative \textsl{h}- and \textsl{v}-covariant derivatives of $X_i$ are given by\\[-10mm]
\begin{equation}
X_{i|\beta}=X_{i|j}B^j_\beta +X_i|_j N^j H_\beta,~~~X_i|_\beta=X_i|_j B^j_\beta.
\end{equation}

A. Rapcs\`ak \cite{ar57} introduced three kinds of hyperplanes. M. Matsumoto \cite{mm85} obtained their characteristic conditions, which are given in the following lemmas\,:
\begin{lemma}
A hypersurface $F^{n-1}$ is a hyperplane of first kind if and only if $H_\alpha=0$ or equivalently $H_0=0$\,.
\end{lemma}
\begin{lemma}
A hypersurface $F^{n-1}$ is a hyperplane of second kind if and only if $H_{\alpha\beta}=0$\,.
\end{lemma}
\begin{lemma}
A hypersurface $F^{n-1}$ is a hyperplane of third kind if and only if $H_{\alpha\beta}=0=M_{\alpha\beta}$\,.
\end{lemma}
\section{The Finsler space $^*\!F^n$=$(M^n,^*\!\!L)$}
~~~~Let us denote $b_iy^i$ by $\beta$, then indicatery property of $h_{ij}$ yield $\dot{\partial}_i\beta=b_i$\,. The quantities corresponding to $^*\!F^n$ is denoted by asterisk over that quantity.\,We shall use following notations  $L_i=\dot\partial_i L=l_i$\,,\,\,\,\,$L_{ij}=\dot\partial_i\dot\partial_j L$\,,\,\,\, $L_{ijk}=\dot\partial_i\dot\partial_j\dot\partial_k L$. From (1.3), we get
\begin{equation}^*L_{ij}=e^\tau (1+\rho-\tau) L_{ij}+\frac{e^\tau}{L}m_im_j\,,\end{equation}
\begin{equation}
\begin{split}
^*L_{ijk}=&\,e^\tau\big(1+\rho-\tau\big)L_{ijk}+\big(\rho-\tau\big)\frac{e^\tau}{L} \big[m_iL_{jk}+ m_jL_{ik}+ m_kL_{ij}\big]\\& -\frac{e^\tau}{L^2}\big[m_jm_kl_i+ m_im_kl_j+m_im_jl_k-m_im_jm_k\big],
\end{split}
\end{equation}
where $\tau=\frac{\beta}{L}$, $m_i=b_i-\tau l_i$\,. The normalised supporting element and the metric tensor of $^*\!F^n$ are obtained as \cite{ma16} \\[-12mm]
\begin{equation}^*l_i=e^\tau\big(m_i+l_i\big),\end{equation}\\[-20mm]
\begin{equation}
^*\!g_{ij}=\nu e^{2\tau}g_{ij} +e^{2\tau}(2\tau^2-\tau-\rho)l_il_j+e^{2\tau}(1-2\tau)(b_il_j+b_jl_i)+2e^{2\tau}b_ib_j \,.
\end{equation}
Differentiating the angular metric tensor $h_{ij}$ with respect to $y^k$, we get
\begin{equation*}
\dot\partial_kh_{ij}=2C_{ijk}-\frac{1}{L}(l_ih_{jk}+l_jh_{ik}),
\end{equation*}\\[-15mm]
which gives\\[-12mm]
\begin{equation}
L_{ijk}=\frac{2}{L}C_{ijk}-\frac{1}{L^2}(h_{ij}l_k+h_{jk}l_i+h_{ki}l_j).
\end{equation}
Using this, the equation (3.2) may be re-written as
\begin{equation}
^*\!C_{ijk}=\nu e^{2\tau}C_{ijk}+\frac{2}{L}e^{2\tau}m_im_jm_k +\frac{1}{2L}e^{2\tau}(2\nu-1)(m_ih_{kj}+m_jh_{ki}+m_kh_{ij})\,,
\end{equation}
where $\nu=1+\rho-\tau\,$.\\
The inverse metric tensor of $^*\!F^n$ is derived as follows\cite{ma16}\,:
\begin{equation}\begin{split}
^*\!g^{ij}=\frac{e^{-2\tau}}{\nu}\Big[g^{ij}-\frac{1}{m^2+\nu}b^ib^j+\frac{\tau-\nu}{m^2+\nu}\big(b^il^j+b^jl^i\big)-\Big\{\frac{\tau-\nu}{m^2+\nu}(m^2+\tau)-\rho\Big\}l_il_j\Big]\,,
\end{split}\end{equation}
where $b$ is magnitude of the vector $b^i=g^{ij}b_j$\,.\\The relation between cartan connection coefficients of $^*\!F^n$ and $F^n$ is given by
\begin{equation} {^*}\!F^i_{jk}=F^i_{jk}+D^i_{jk}\,.\end{equation}\\[-10mm]
The expressions for $D^i_{00}$, $D^i_{0k}$ and $D^i_{jk}$ are given by \cite{ma16}
\begin{equation}\begin{split}
D^i_{00}=&\frac{L}{\nu e^\tau}\Big[\frac{e^\tau}{L}\beta_{|0}m^i+2e^\tau F^i_0\Big]
+l^i\Big[E_{00}-\frac{L}{e^\tau }\big(m^2+\nu \big)^{-1}\big(\frac{e^\tau}{L}\beta_{|0}m^2+2e^\tau F_{\beta 0}\big)\Big]\\& -\frac{m^iL}{\nu e^\tau}\big(m^2+\nu \big)^{-1}\Big[\frac{e^\tau}{L}\beta_{|0}m^2+2e^\tau F_{\beta 0}\Big]\,,
\end{split}\end{equation}
\begin{equation}\begin{split}
D^i_{0j}=\frac{LG^i_j}{\nu e^\tau}+\frac{l^i}{e^\tau}\Big[G_j\,-L\big(m^2+\nu \big)^{-1} G_{\beta j}]-\frac{m^iL}{\nu e^\tau}\big(m^2+\nu \big)^{-1} G_{\beta j}\,,
\end{split}\end{equation}
\begin{equation}\begin{split}
D^j_{ik}=\frac{LH^j_{ik}}{\nu e^\tau}+ \frac{l^j}{e^\tau}\Big[{H_{ik}}-{L\big(m^2+\nu \big)^{-1}H_{\beta ik}}\Big]\,-\frac{m^jL}{\nu e^\tau}\big(m^2+\nu \big)^{-1} H_{\beta ik} \,,
\end{split}\end{equation}
where 
\begin{equation}\begin{split}
 2G_{ij}=&\frac{e^\tau}{L}\big(\beta_{|j}m_i-\beta_{|i}m_j\big)+2e^\tau F_{ij}-\nu e^\tau L_{ijr}D^r_{00}-A_{ijk}y^k-\frac {e^\tau}{L^2}m_rm_im_jD^r_{00}\\&+\frac{(\nu -1)}{L}e^\tau \beta_{|0}L_{ij}+\frac{e^\tau}{L^2}B_0m_im_j+e^\tau \rho_0 L_{ij}\,,
\end{split}\end{equation}
\begin{equation}
G_j=e^\tau\big(E_{j0}-F_{j0}\big),
\end{equation}
\begin{equation}\begin{split}
2H_{jik}=& -\nu e^\tau\Big[L_{ij r}D^r_{0k}+L_{jkr}D^r_{0i}-L_{kir}D^r_{0j}\Big]+A_{jki} +A_{kij} -A_{ijk} \\&
-\frac{e^\tau}{L^2}\Big[m_im_jm_rD^r_{0k}+m_jm_km_rD^r_{0i}-m_km_im_rD^r_{0j}\Big]\\&+(\nu-1)\frac{e^\tau}{L}\big(\beta_{|k}L_{ij}+\beta_{|i}L_{jk}-\beta_{|j}L_{ki}\big)+e^\tau  \Big[\rho_kL_{ij}+\rho_iL_{jk}-\rho_j L_{ki}\Big]\\&
+\frac{e^\tau}{L^2}\Big[ \beta_{|k}m_im_j+\beta_{|i}m_jm_k-\beta_{|j}m_km_i\Big]\,, \end{split}\end{equation}
\begin{equation}\begin{split}
2H_{ik}=&\frac{e^\tau}{L}\big(\beta_{|k}m_i+\beta_{|i}m_k\big)+e^\tau E_{ik}-\Big[\nu e^\tau L_{ir}+\frac{e^\tau}{L}m_im_r\Big]D^r_{0k}\\&-\Big[\,\nu e^\tau L_{kr}+\frac{e^\tau}{L}m_km_r\Big]D^r_{0i}\,,
\end{split}\end{equation}
\begin{equation}
A_{ijk}=\frac{e^\tau}{L}D^r_{0k}\mathfrak{S}_{(rij)}\Big[(\nu-1) m_rL_{ij}-\frac{m_im_jl_r}{L}\Big]\,.
\end{equation}
Here we have used $m_im^i=m^2=m^ib_i$ and $H^j_{ik}=g^{jm}H_{mik}$.\,Also we note that $E_{00}=E_{ij}y^iy^j=b_{i|j}y^iy^j=(b_iy^i)_{|j}y^j=\beta_{|0}$,$~$$F^i_0=g^{ij}F_{j0}$\,. The subscript `0' denotes the contraction by supporting element $y^i$, unless otherwise stated. The subscript `$\beta$' denotes the contraction by supporting element $b^i$ .\,$\mathfrak{S}_{(ijk)}$ denote cyclic interchange of indices $i$,\,$j$,\,$k$ and summation.\\
Now, we state a Lemma which are used later.
\begin{lemma} \textbf{$\cite{mp15}$}
If the $h$-vector $b_i$ is gradient then the scaler $\rho$ is constant\,.
\end{lemma}
\section{The Hypersurface $^*\!F^{n-1}$ of the space $^*\!F^{n}$}
~~~~Let us consider  Finslerian hypersurfaces $F^{n-1}=(M^{n-1},\underline L(u,v))$ of $F^n$ and $^*\!F^{n-1}=(M^{n-1},\underline {^*\!L}(u,v))$ of\, $^*\!F^n$. Let $N^i$ be the unit normal vector at a point of $F^{n-\!\!1}$.\,The functions $B^i_\alpha(u)$ may be considered as component of $(n-1)$ linearly independent vectors tangent to $F^{n-1}$ and they are invariant under h-exponential change of Finsler metric.\,The unit normal vector $^*\!N^i(u,v)$ of $^*\!F^{n-1}$ is uniquely determined by
\begin{equation}
{^*\!g}_{ij}B^i_\alpha {^*\!N^j}=0\,, ~~~~~~~ ^*\!g_{ij}{^*\!N^i}{^*\!N^j}=1\,.
\end{equation}
The inverse projection factors $^*\!B^\alpha_i(u,v)$ of $B^i_\alpha$ along $^*\!F^{n-1}$ are defined as\\[-10mm] 
\begin{equation}
{^*\!B}^\alpha_i={^*\!g}^{\alpha\beta}{^*\!g}_{ij}B^j_\alpha\,,
\end{equation}\\[-10mm]
where $^*\!g^{\alpha\beta}$ is the inverse of metric tensor $^*\!g_{\alpha\beta}$ of $^*\!F^{n-1}$.\\
From (4.2), it follows that
\begin{equation}
B^i_\alpha {^*\!B}^\beta_i=\delta^\beta_\alpha\,, ~~~~B^i_\alpha {^*\!N_i}=0 \,,~~~~~{^*\!N}_iB^i_\alpha=0\,,~~~~~ {^*\!N_i}{^*\! N^i}=1\,,
\end{equation}\\[-10mm]
and further\\[-10mm]
\begin{equation}
B^i_\alpha {^*\!B}^\beta_j+ {^*\!N^i}{^*\! N_j}=\delta^i_j\,.
 \end{equation}

Now, Transvection of (2.3) by $v^\alpha$  gives
\begin{equation}
y_jN^j=0\,.
\end{equation}
Transvecting (3.4) by $N^iN^j$ and by using (4.5), we have
\begin{equation}
^*\!g_{ij}N^iN^j=\nu e^{2\tau} +2e^{2\tau}(b_iN^i)^2 ,
\end{equation}\\[-10mm]
this implies that\\[-10mm]
\begin{equation}
\frac{N^j}{{e^{\tau}}\sqrt{\nu +2(b_iN^i)^2}}
\end{equation}\\[-10mm]
is unit vector.\\
Also, Transvecting (3.4) by $B^i_\alpha N^j$ and using (4.5), gives us
\begin{equation}
^*\!g_{ij}B^i_\alpha N^j=(b_jN^j)e^{2\tau}\Big\{(1-2\tau)l_iB^i_\alpha +2b_iB^i_\alpha\Big\} \,.
\end{equation}
This shows that $N^j$ is normal if and only if R.H.S. of equation (4.8) is zero. Since \,\,$e^{2\tau}\Big\{(1-2\tau)l_iB^i_\alpha +2b_iB^i_\alpha\Big\}$\,\, can not be zero, otherwise transvection of\,\, $e^{2\tau}\Big\{(1-2\tau)l_iB^i_\alpha +2b_iB^i_\alpha\Big\}$ by $v^\alpha$\,\, gives $L=0$, which is not possible. Hence $N^j$ is normal to $^*\!F^{n-1}$ if and only if $b_jN^j=0$.\\From (4.7) and (4.8), we may state that
\begin{equation}
{^*\!N^i}=\frac{N^i}{{e^\tau}\sqrt \nu}
\end{equation}\\[-10mm]
is unit normal vector of $^*\!F^{n-1}$.\\Which in view of (3.4) and (4.5), gives\\[-10mm]
\begin{equation}
{^*\!N_i}={N_i}{{e^\tau}\sqrt \nu}\,.
\end{equation}\\[-10mm]
Thus, we have\,:
\begin{theorem}
Let $^*\!F^n$ be the Finsler space obtained from $F^n$ by $\textsl{h}$-exponential change given by (1.3). Further if $^*\!F^{n-1}$ and $F^{n-1}$ are the hypersurfaces of these spaces.\,Then the vector $b_i$ is tangential to hypersurface $F^{n-1}$ if and only if every vector normal to $F^{n-1}$ is also normal to $^*\!F^{n-1}$. And then the normal vector is given by (4.9).
\end{theorem}
\indent Let $b_i$ is gradient vector, i.e. $b_{j|i}=b_{i|j}$, then\\[-10mm]
\begin{equation}
F_{ij}=0\,,
\end{equation}\\[-10mm]
which in view of Lemma (3.1), gives\\[-10mm]
\begin{equation}
\rho_i=0\,.
\end{equation}\\[-10mm] 
Now, if $b_i$ is tangent to hyperplane $F^{n-1}$ $i.\,e.\,$
\begin{equation}
b_jN^j=0\,.
\end{equation}\\[-10mm]
Using (4.5), (4.11) and (4.13), we have
\begin{equation}
D^i_{00}N_i=0\,.
\end{equation}
The normal curvature tensor $^*\!H_\alpha$ for hypersurface $^*\!F^{n-1}$ is given by
\begin{equation*}
{^*\!}H_\alpha={^*\!}N_i(B^i_{0\alpha}+{^*\!}G^i_j B^j_\alpha)\,,
\end{equation*}
by use of (2.8) and (4.9), above equation becomes
\begin{equation}
{^*\!}H_\alpha={\sqrt \nu}\,\, e^{\tau}\Big(H_\alpha+N_iD^i_{0j}B^j_\alpha\Big)\,,
\end{equation}
 which on transvection by $v^\alpha$ and using $(4.14)$, gives
\begin{equation}
{^*\!}H_0={\sqrt \nu}\,\,e^{\tau}H_{0}\,.
\end{equation}
\noindent Thus in view of Lemma (2.1), we have\,:
\begin{theorem}
Let the $h$-vector $b_i$ be a gradient and tangent to hypersurface $F^{n-1}$.\,Then the hypersurface $F^{n-1}$ is a hyperplane of first kind if and only if hypersurface $^*\!F^{n-1}$ is hyperplane of first kind.
\end{theorem}
Taking the relative $\textsl{h}$-covariant differentiation of (4.13) with respect to the Cartan connection of $F^{n-1}$, we get
\begin{equation*}
b_{i|\beta}N^i + b_i N^i_{\beta}=0\,.
\end{equation*}
Using (2.11) and (2.12), the above equation gives
 \begin{equation*}
(b_{i|j}B^j_\beta + b_i|_j N^j H_\beta)N^i - b_iH_{\alpha \beta} B^\alpha_j g^{ij}=0\,.
\end{equation*}
Travecting by $v^\beta$ and using (2.9), we get
\begin{equation*}
b_{i|0}N^i=(H_\alpha + M_\alpha H_0)B^\alpha_j b^j - b_i|_j H_0 N^i N^j\,.
\end{equation*}
For the hypersurface to be first kind, $H_0=0=H_\alpha$.\,Then above equation reduces to $b_{i|0} N^i=0$\,. If the vector $b_i$ is gradiant,\,\textit{i.e.} $b_{i|j}=b_{j|i}$, then we get 
\begin{equation*}
E_{i0}N^i=b_{i|0}N^i =\beta_i N^i\,.
\end{equation*} 
The tensors $D^i_{00},\,D^i_{oj},\,G_{ij}$ and $G_j$ satisfies the following, which can be easily verified\,:
\begin{equation}\begin{split}
& D^i_{00}N_i=0\,,~~~~D^r_{0j}L_{jr}N^j=0\\
&L_{ijr}D^r_{00}=\big(E_{00}-(m^2+\nu)^{-1}\beta_0m^2\big)\Big[\big(\frac{2\rho}{L^2\nu}-\frac{1}{L^2}\big)h_{ij}-\frac{1}{L^2 \nu}(m_jl_i+m_il_j)\Big]\,,\\
&G_{ij}N^iB^j_\alpha=0\,,~~~~
D^i_{0j}N_iB^j_\alpha=0\,,~~~~G_jN^j=0\,,~~~~G_{ij}b^iN^j=0\,,~~~~\\&D^i_{0j}b_iN^j=0\,,~~~~D^i_{0j}N^jB^k_\alpha h_{ik}=0\,,~~~~D^r_{0j}l_rN^j=0\,,~~~~G^r_j l_rN^j=0\,.
\end{split}\end{equation}
The second fundamental \textsl{h}- tensor $^*\!H_{\alpha\beta}$ for hyperplane $^*\!F^{n-1}$ is given by
\begin{equation*}
{^*\!H}_{\alpha\beta}-{^*\!M}_\alpha{^*\!H}_\beta={^*\!}N_i(B^i_{\alpha \beta}+{^*\!}F^i_{jk}B^j_\alpha B^k_\beta)\,, 
\end{equation*}
then by use of (2.7),\,(3.8) and (4.10), above equation gives
\begin{equation}
{^*\!H}_{\alpha\beta}-{^*\!M}_\alpha{^*\!H}_\beta={{e^{\tau}}\sqrt \nu}\Big[H_{\alpha\beta}+N_iD^i_{jk}B^j_{\alpha}B^k_{\beta}\Big] -{e^{\tau}\sqrt\nu}M_\alpha H_\beta \,.
\end{equation}
Contracting  (3.11) by $B^i_\alpha B^k_\beta N_j$ and using $m^j N_j=0$,\,\,$l^jN_j=0$, we get 
\begin{equation*}\begin{split}
D^j_{ik}B^i_{\alpha}B^k_{\beta}N_j=\frac{L}{\nu{e^{\tau}}}H^j_{ik}B^i_{\alpha}B^k_{\beta}N_j=-\frac{L}{2\nu e^\tau}H_{jik}N^jB^i_{\alpha}B^k_{\beta}\,,
\end{split}\end{equation*}
which in view of (3.14) and (4.17), gives
\begin{equation}\begin{split}
D^j_{ik}B^i_{\alpha}B^k_{\beta}N_j=-\frac{L}{2} \Big[L_{ij r}D^r_{0k}+L_{jkr}D^r_{0i}-L_{kir}D^r_{0j}\Big]N^jB^i_{\alpha}B^k_{\beta}\,.
\end{split}\end{equation}

Now we calculate each terms of the above equation separately.\\
Transvecting (3.12) by $N^j$, we have
\begin{equation}
  G_{ij}N^j= \mu N_i\,,\end{equation}\\[-12mm]
where\\[-10mm]
\begin{equation*}\begin{split}
\mu=\frac{1}{2L^2}\Big[-{e^\tau}\big(E_{00}-(m^2+\nu)^{-1}\beta_0m^2\big)\big(2\rho-\nu\big)-{e^\tau}\nu D^r_{00}m_r+\nu({\nu-1})e^\tau \beta_0\Big]\,.
\end{split}\end{equation*}
Contracting $L_{ijr}$ by $N^jB^i_{\alpha}B^k_{\beta} D^r_{0k}$ and using (1.2)\,,(3.10) and above equation, we obtain
\begin{equation}
L_{kir}N^jB^i_{\alpha}b^k_{\beta}D^r_{0j}=\frac{2\mu}{{ \nu e^{\tau}}}M_{\alpha\beta}\,,
\end{equation} 
Transvecting $L_{ijr}$ by $N^jB^i_{\alpha}B^k_{\beta}$ and using (1.2),\,(3.10) and (3.12), we get 
\begin{equation}\begin{split}
L_{ij r}N^jB^i_{\alpha}B^k_{\beta} D^r_{0k}=\frac{2}{{\nu e^{\tau}}}\Big[\lambda M_{\alpha \beta}-\frac{e^\tau}{2L}\beta_r C^r_{ij}N^jB^i_\alpha B^k_\beta m_k\Big]\,,
\end{split}\end{equation}\\[-10mm]
where\\[-10mm]
\begin{equation*}
\lambda=\frac{1}{2L^2}\Big[-e^\tau\big(E_{00}-(m^2+\nu)^{-1}\beta_{0}m^2\big)\Big(2\rho-\nu\Big)-{e^\tau}\nu\,(\nu-1) D^s_{00}m_s+\nu\,(\nu-1)e^\tau \beta_{0}\Big]\,.
\end{equation*}
Similarly, transvecting $L_{kjr}$ by $N^jB^i_{\alpha}B^k_{\beta} D^r_{0i}$ and using $M_{\alpha \beta}= M_{\beta \alpha}$, we have 
\begin{equation}
L_{kj r}N^jB^i_{\alpha}B^k_{\beta}D^r_{0i}=\frac{2}{\nu e^{\tau}}\Big[ \lambda M_{\alpha \beta}-\frac{e^\tau}{2L}\beta_r C^r_{ij}N^jB^i_\beta B^k_\alpha m_k\Big]\,.
\end{equation}
Plugging (4.21),\,(4.22),\,(4.23) in equation (4.19), we obtain
\begin{equation}
D^j_{ik}N_jB^i_{\alpha}B^k_{\beta}= \frac{L(\mu-2\lambda)}{e^\tau\,\nu}M_{\alpha \beta}+\frac{e^\tau}{2L}\beta_r C^r_{ij}\Big[N^jB^i_\alpha B^k_\beta m_k+N^jB^i_\beta B^k_\alpha m_k \Big]\,.
\end{equation}
Now, suppose that $h$-vector $b_i$ satisfies the condition 
\begin{equation}
b_{r|0}\,C^r_{ij}=\kappa\,h_{ij}\,,
\end{equation}\\[-10mm]
then\\[-8mm]
\begin{equation}\beta_r\, C^r_{ij}=\kappa\,h_{ij}\,,
\end{equation}
where $\kappa$ is a scalar function.\\
So, using $h_{ij}B^i_\alpha N^j=0$,  equation (4.24) yields
\begin{equation}
D^j_{ik}B^i_{\alpha}B^k_{\beta}N_j=\frac{L(\mu-2\lambda)M_{\alpha \beta}}{\nu e^{\tau}}\,.
\end{equation}
And then (4.18) becomes
\begin{equation}
{^*\!H}_{\alpha \beta}-{^*\!M}_\alpha {^*\!H}_{\beta}=e^\tau\sqrt{\nu}\Big[H_{\alpha \beta}+\frac{L(\mu-2\lambda)}{\nu e^\tau}\,M_{\alpha \beta}\Big]-e^{\tau}\sqrt\nu\,M_\alpha H_\beta\,\,. 
\end{equation}
Next, transvecting  (3.6) by $B^i_\alpha B^j_\beta N^k$ and using (2.10), we have
\begin{equation}
{^*\!}M_{\alpha\beta}= {\sqrt \nu}\,\,e^{\tau}M_{\alpha\beta}\,.
\end{equation}
\noindent Thus from (2.28) and (2.29), we have\,:
\begin{theorem}
For the exponential change with an $h$-vector, let the $h$-vector $b_i$ be a gradient and tangential to hypersurface $F^{n-1}$ and satisfies condition (4.25).\,Then
\begin{enumerate}
\item {${^*\!F^{n-1}}$ is a hyperplane of second kind if $F^{n-1}$ is hyperplane of second kind and $M_{\alpha\beta}=0$.}
\item{${^*\!F^{n-1}}$ is a hyperplane of third kind if $F^{n-1}$ is hyperplane of third kind.}
\end{enumerate}
\end{theorem}
\section{Example}
~~~~A Finsler space $F^n$ is called $^*\!P$-Finsler space  if the $(v)hv$-torsion tensor $P^r_{ij}$ satisfies\\[-12mm]
\begin{equation}
P^r_{ij}:=C^r_{ij|0}= \lambda \,C^r_{ij}\,.
\end{equation}\\[-12mm]
Taking $h$-covariant derivative of (1.1) and using $L_{|k}=0=h_{ij|k}$ and $\rho_i=0$, we get
\begin{equation}
b_{r|k}C^r_{ij}+b_r C^r_{ij|k}=0\,.
\end{equation}\\[-8mm]
Contracting the above equation by $y^k$ and using (5.1), we get\\[-8mm]
\begin{equation*}
b_{r|0}\,C^r_{ij}+\lambda\, b_r\, C^r_{ij}=0\,,
\end{equation*}\\[-8mm]
which in view of (1.1),\textsc{} becomes
\begin{equation}
b_{r|0}\,C^r_{ij}=\kappa\,h_{ij}\,,~~ \kappa=-\frac{\lambda \rho}{L},
\end{equation}\\[-12mm]
which is required condition (4.25). Thus, we have\,:
\begin{theorem}
For the exponential change with an $h$-vector, let the $h$-vector $b_i$ be a gradient and tangential to hypersurface $F^{n-1}$ of a $^*\!P$-Finsler space $F^n$.\,Then
\begin{enumerate}
\item {${^*\!F^{n-1}}$ is a hyperplane of second kind if $F^{n-1}$ is hyperplane of second kind and $M_{\alpha\beta}=0$}.
\item{${^*\!F^{n-1}}$ is a hyperplane of third kind if $F^{n-1}$ is hyperplane of third kind.}
\end{enumerate}
\end{theorem}

A Landsberg space is $^*\!P$-Finsler space for $\kappa =0$.\\ 
Thus, we have\,:
\begin{corollary}
For the exponential change with an $h$-vector, let the $h$-vector $b_i$ be a gradient and tangential to hypersurface $F^{n-1}$ of a Landsberg space $F^n$.\,Then
\begin{enumerate}
\item {${^*\!F^{n-1}}$ is a hyperplane of second kind if $F^{n-1}$ is hyperplane of second kind and $M_{\alpha\beta}=0$}.
\item{${^*\!F^{n-1}}$ is a hyperplane of third kind if $F^{n-1}$ is hyperplane of third kind.}
\end{enumerate}

\end{corollary}
\section*{Discussion}
~~~~Gupta and Pandey \cite{mp16} have proved that for Kropina change with an $h$-vector (let the $h$-vector $b_i$ be a gradient and tangential to hypersurface $F^{n-1}$ and satisfies condition $\beta_r C^r_{ij}=0$)\,,\\[-10mm]
\begin{center}
\textbf{{${^*\!F^{n-1}}$ is a hyperplane of third kind if $F^{n-1}$ is hyperplane of third kind.}} 
\end{center}

In present paper, authors proved that for exponential change with an $h$-vector (same conditions)\,,\\[-12mm]
\begin{center}
\textbf{{${^*\!F^{n-1}}$ is a hyperplane of third kind if $F^{n-1}$ is hyperplane of third kind.}
} 
\end{center}

Notice that Kropina change with an $h$-vector is finite in nature (in the sense that number of terms) whereas exponential change with an $h$-vector is infinite in nature, although in both cases (finite and infinite) same result holds.

\textbf{The question is that}
\textit{Is there any particular type of change with an $h$-vector (same conditions) for which
${^*\!F^{n-1}}$ is a hyperplane of third kind if $F^{n-1}$ is hyperplane of third kind\,?
}


\small\begin{thebibliography}{}
\bibitem [1]{ec34} E.\,Cartan,\, \textit{Les espaces de Finsler}.Actualits 79, specialement XI.,\,1934.
\bibitem[2]{ma16} M.$\,$K.$\,$Gupta and Anil K.\,Gupta,\,$\textsl{h}$-exonential change of Finsler metric, arXiv\,:\,1603. 05434v1 [math.\,DG]\,17 Mar 2016.
\bibitem[3]{mp16} M.$\,$K.$\,$Gupta and P.$\,$N.$\,$Pandey,\,\,Hypersurface of a Finsler space subjected to  a Kropina change with an $h$-vector,\,\,arXiv\,:\,1502.\,0186v1[math.\,DG]\,\,6 Feb 2015.
\bibitem[4]{mp08} M.$\,$K.$\,$Gupta and P.$\,$N.$\,$Pandey,\, On hyperspaces of a Finsler space with a special metric,\,\textit{Acta Math. Hunger.},\,120(1-2),\,(2008),\,165-177.

\bibitem[5]{mp09} M.$\,$K.$\,$Gupta and P.$\,$N.$\,$Pandey,\, Hyperspaces of conformally and \textsl{h}-conformally related Finsler spaces,\,\textit{Acta Math. Hunger.},\,123(3),\,(2009),\,257-264.

\bibitem[6]{mp15} M.$\,$K.$\,$Gupta and P.$\,$N.$\,$Pandey,\, Finsler space subjected to a Kropina change with an \textsl{h}-vector,\,\textit {Facta Uni. Series:\,Maths. and Infor.},\,30(4),\,(2015),\,513-525.

\bibitem[7]{hi80} H.$\,$Izumi\,, Conformal transformation of Finsler spaces II. An \textsl{h}-conformally flat Finsler space,\,\textit {Tensor N.S.},\, 34,\,(1980),\, 337-359.
\bibitem[8]{mm85} M.$\,$ Matsumoto\,.The induced and intrinsic Finsler connections of a hypersurface and Finslerian projective geometry.\,J.\,Math\,Kyoto Univ.,\,\,25(3),\,107-144,\,1985.

\bibitem[9]{bn90} B.$\,$N.$\,$Prasad,\, On the torsion tensors $R_{hjk}$ and $P_{hjk}$ of Finsler spaces with a metric $ds= (g_{ij}(dx)dx^idx^j)^\frac{1}{2}+b_i(x,y)dx^i$, \textit{Indian J.\,Pure Appl.\, Math.},\,\,21(1),\,1990,\,27-39.


\bibitem[10]{hb12} B.$\,$N.\,Prasad,\,H.$\,$S.\,Shukla and O.$\,$P.\,Pandey\,, Exponential Change of Finsler metric, \textit {Int.\,J.Contemp.\,Math.\,Sciences},\,7,\,(2012),\,2253-2263.
\bibitem[11]{ar57} A.\, Rapcs\`ak,\,Eine neue Charakterisierung Finslerscher Räume Skalarer und konstanter Krümmung und projektiv-ebene Räume\,,\textit{Acta Math. Acad. Sci. Hungar}\,,8(3),\,(1957),\,1-8.
\bibitem[12]{yy06} YU\, Yao-yong,\, YOU \,Ying,\,Projectively flat exponential Finsler metric,\,\textit{Journal of Zhejiang University Sciencea}
\,7(6),\,(2006),\,1068-1076.
\end{thebibliography}
\end{document}